\documentclass{amsart}

\usepackage{amssymb}
\usepackage{algorithm}
\usepackage{algpseudocode}
\usepackage{url}
\pdfoutput=1

\theoremstyle{plain}

\theoremstyle{definition}

\theoremstyle{remark}
 
 \numberwithin{equation}{section}

\providecommand{\abs}[1]{\left\lvert#1\right\rvert}
\providecommand{\floor}[1]{\left\lfloor#1\right\rfloor}

\renewcommand{\leq}{\leqslant}

\title[Searching for a counterexample to Kurepa's Conjecture]{SEARCHING FOR A COUNTEREXAMPLE TO KUREPA'S CONJECTURE}



\author[Andreji\'c]{\bfseries Vladica Andreji\'c}

\address{
Faculty of Mathematics \\
University of Belgrade  \\
Belgrade\\
Serbia} \email{andrew@matf.bg.ac.rs}

\author[Tatarevic]{\bfseries Milos Tatarevic}
\address{Alameda, CA 94501}
\email{milos.tatarevic@gmail.com}

\thanks{This work is partially supported by the Serbian Ministry of Education and Science, project No. 174012}

\begin{document}

\vspace{18mm}
\setcounter{page}{1}
\thispagestyle{empty}

\begin{abstract}
Kurepa's conjecture states that there is no odd prime $p$ that
divides $!p=0!+1!+\cdots+(p-1)!$. We search for a counterexample to
this conjecture for all $p<2^{34}$. We introduce new optimization
techniques and perform the computation using graphics processing
units. Additionally, we consider the generalized Kurepa's left
factorial given by $!^{k}n=(0!)^k +(1!)^k +\cdots+((n-1)!)^{k}$, and show
that for all integers $1<k<100$ there exists an odd prime $p$ such
that $p\mid !^k p$.
\end{abstract}

\maketitle

\section{Introduction}

\DJ uro Kurepa defines an arithmetic function !n, known as the left factorial,
by $!n=0!+1!+\cdots+(n-1)!$ \cite{K1}. He
conjectures that $\text{GCD}(!n,n!)=2$ holds for all integers $n>1$. This conjecture,
originally introduced in 1971 \cite{K1}, is also listed in \cite[Section B44]{Gu}
and remains an open problem. For additional information, the reader can consult the
expository article \cite{IM}.

It is easy to see that the statement above is equivalent to the
statement that $!n$ is not divisible by $n>2$, which can be
reduced to primes. Thus, Kurepa's conjecture states that there is
no odd prime $p$ such that $p$ divides $!p$.

There have been numerous attempts to solve this problem, mostly by
searching for a counterexample. The results of these attempts are listed in
the accompanying table.

\vspace{5mm}
\begin{center}
\begin{tabular}{|l|l|c|}
  \hline
  Upper Bound & Author & Year\\
  \hline
  $p<3\cdot10^5$ & \v Z. Mijajlovi\'c \cite{Mi} & $1990$\\
  $p<10^6$ & G. Gogi\'c \cite{Go} & $1991$\\
  $p<3\cdot10^6$ & B. Male\v sevi\'c \cite{Ma} & $1998$\\
  $p<2^{23}$ & M. \v Zivkovi\'c \cite{Z} & $1999$\\
  $p<2^{26}$ & Y. Gallot \cite{Ga} & $2000$\\
  $p<1.44\cdot 10^{8}$ & P. Jobling \cite{J} & $2004$\\
 \hline
\end{tabular}
\end{center}
\vspace{5mm}

In 2004, Barsky and Benzaghou published a proof of Kurepa's
conjecture \cite{BB1}, but owing to irreparable calculation errors
the proof was retracted \cite{BB2}. Belief that the problem was
solved may explain why there have been no recent attempts
to find a counterexample.

In this article, we extend the search to all primes $p<2^{34}$.
The number of arithmetic operations required to search for a
counterexample in this interval is greater by a factor of approximately
$11300$ compared to the latest published results.
To achieve this, we improve the algorithm by
reducing the number of required modular reductions by a factor
of $2$ and perform our calculations on graphics processing units
(GPUs).

We have calculated and recorded the residues $r_{p}\!\!=\,\,!p \bmod p$ for
all primes $p<2^{34}$. We confirm all $r_p$ given in results from
previous searches and find no counterexample ($r_{p}=0$). It
is worth noting that $r_p=0$ is not the only interesting residue.
For example, $r_p=1$ is also mentioned in \cite[Section B44]{Gu},
with two solutions $p=3$ and $p=11$, as well as the remark that
there are no other solutions for $p<10^6$. Our search verifies
that there are no new solutions to $p\mid 1!+2!+\cdots+(p-1)!$
for $p<2^{34}$. However, we managed to find a new solution to
$r_{p}=2$, which means that $p\mid 2!+\cdots+(p-1)!$ for
$p=6855730873$. This is the first such prime after $p=31$ and
$p=373$.

Finally, for a positive integer $k$, we consider the generalized Kurepa's left factorial
$$!^k n=(0!)^k+ (1!)^k+\cdots+((n-1)!)^k,$$
where $!^1 k=!k$. Similarly, we can ask
whether some odd prime $p$ divides $!^k p$. Such considerations
are known from Brown's MathPages \cite{Br}, where the analog
of Kurepa's conjecture for $k=5$ is introduced as an open problem.
We find a counterexample, $p\mid !^5 p$ for $p=9632267$,
and present further examples where $p\mid !^k p$ for all $1<k<100$.

\section{Algorithmic considerations}

Wilson's theorem states that $(p-1)!\equiv -1\pmod p$ for all primes $p$.
It is easy to show that
\begin{equation}\label{wilson}
(p-k)!(k-1)!\equiv(-1)^k \pmod p,
\end{equation}
for all primes $p$ and all $1\leq k\leq p$. Using (\ref{wilson}),
we can obtain
\begin{equation}\label{levi}
!p\equiv\sum_{i=0}^{p-1}(-1)^{i}\frac{(p-1)!}{i!}\pmod p,
\end{equation}
or, using the theory of derangements,
$$!p\equiv\floor{\frac{(p-1)!}{e}} +1\pmod p.$$

The value $r_p$ can be calculated using recurrent formulas.
Regroup the parentheses in the definition of $!p$ to get
$$!p=1+1(1+2(1+3(1+\cdots(1+(p-2)(1+(p-1)))\cdots))).$$
This formula shows that if we define $A_{i}$ by $A_{1}=0$ and
$A_{i}=(1-iA_{i-1})\bmod p$ for $i>0$, then $A_{p-1}=r_{p}$.
Likewise, regrouping (\ref{levi}) in the same way shows that if we
define $B_{i}$ by $B_{1}=0$ and $B_{i}=((-1)^i+iB_{i-1})\bmod p $
for $i>1$, then $B_{p-1}=r_{p}$.  Both formulas appear in
\cite{IM}.  It is worth mentioning that the sequence
generated by the same recurrence as $B_{i}$ but without modular
reduction is a well-known combinatorial sequence representing
the number of derangements \cite{OEIS}.
Similarly, we can use the definition $D_{1}=C_{1}=1$,
$D_{i}=iD_{i-1}$, $C_{i}=(C_{i-1}+D_{i-1})\bmod p$ and get
$r_{p}=C_{p}$, as proved in \cite{Ga}.

The best-known algorithm for computing a single value $r_{p}$ has
time complexity $O(p)$, which implies that the time required to search
for a counterexample for all $p<n$ is $O(n^{2}/\ln n)$.
In practice, computer hardware performs integer arithmetic with a fixed precision,
so the product of two operands often needs to be reduced modulo $p$.
All previous attempts to resolve $r_{p}$ required
$p$ modular reductions and a few additional multiplications and additions.

Modern hardware is optimized to perform fast multiplication and
addition on integers and floats, while the biggest concern is
division which is, on average, orders of magnitude slower. Algorithms
that compute $a\cdot b\pmod p$ quickly, such as Montgomery reduction
\cite{Mo} or inverse multiplication \cite{Ba}, use multiplication
and addition to replace division, resulting in execution being several times
faster.

Since we tested $p < 2^{34}$, we could not use the
32-bit integer arithmetic. Instead, we used double precision floating point
hardware in a novel way to compute with integers up to $h = 2^{51}$.
Note that the IEEE 754 Standard for double precision floating point
arithmetic specifies a 52-bit significand.

Our algorithm uses a modified approach that allows us to reduce the
number of required modular reductions by a factor of $2$.
While $p<h$, we can avoid modular reduction in cases where we need to
perform $p+a$ for any $a<h-p$, or when we multiply $p$ by a
constant $a<\frac{h}{p}$.

Let us regroup $!p$ in the following way:
\begin{equation*}
!p=\sum_{i=0}^{p-1}i!=1+(1!+2!)+(3!+4!)+\cdots+((p-4)!+(p-3)!)+((p-2)!+(p-1)!).
\end{equation*}
Because $(p-1)!\equiv-1\pmod p$ and $(p-2)!\equiv1\pmod p$, we have
\begin{align*}
!p&\equiv1+\sum_{i=1}^{\frac{p-3}{2}}\left((p-2i-2)!+(p-2i-1)!\right)
\pmod
p\\
&\equiv1+\sum_{i=1}^{\frac{p-3}{2}}(p-2i-2)!(1+(p-2i-1)) \pmod p\\
&\equiv1-2\sum_{i=1}^{\frac{p-3}{2}}i(p-2i-2)! \pmod p.
\end{align*}
From (\ref{wilson}) we have $(p-2i-2)!(2i+1)!\equiv(p-2)!\pmod p$, and
therefore
\begin{equation}\label{evoga}
r_p\equiv 1-2\sum_{i=1}^{\frac{p-3}{2}}i\frac{(p-2)!}{(2i+1)!}
\pmod p.
\end{equation}
The sequence
$$s_k=\sum_{i=1}^{k}\frac{i}{k}\frac{(2k+1)!}{(2i+1)!}$$
can be written recursively as
$$s_{i}=(2i-2)(2i+1)s_{i-1}+1$$
with $s_{1}=1$, because of the induction step
\begin{align*}
(2k-2)(2k+1)s_{k-1}+1&=
2(k-1)(2k+1)\sum_{i=1}^{k-1}\frac{i}{k-1}\frac{(2k-1)!}{(2i+1)!}+1\\
&= \sum_{i=1}^{k-1}\frac{i}{k}\frac{(2k+1)!}{(2i+1)!}+1=s_k.
\end{align*}
Thus (\ref{evoga}) gives
$$r_p\equiv 1-2 \frac{p-3}{2}\sum_{i=1}^{\frac{p-3}{2}}\frac{i}{\frac{p-3}{2}}\frac{(2\frac{p-3}{2}+1)!}{(2i+1)!}
\equiv 1+3s_{\frac{p-3}{2}} \pmod p.$$ We provide the final form
of $s_{i}$ for $1\leq i\leq \frac{p-3}{2}$ using two new
sequences:
\begin{equation}\label{s_i}
s_{i}=(m_{i-1}s_{i-1}+1)\bmod p,\qquad m_{i}=m_{i-1}+k_{i-1},\quad
k_{i}=k_{i-1}+8,
\end{equation}
with initial values $s_{1}=1, m_{1}=10, k_{1}=18$, which gives
$r_p\equiv 1+3s_{\frac{p-3}{2}} \pmod p$.

To ensure that $m_{i}$ remains below $h$, we do not need to perform modular reduction until after
a large number of iterations. If we start the loop with
$k_{i}\leq p$, then while $k_{i}$ increases linearly by $8$, $m_{i}$
will exceed $h$ after approximately $\frac{h}{p}$
iterations. Because we examine primes less than $2^{34}$, we can perform
approximately $2^{17}$ iterations before a modular reduction.

Although we were able to reduce the number of operations required to
calculate $r_{p}$, the time complexity of our algorithm is still $O(n^{2}/\ln n)$.

\section{Machine considerations}

In the last few years the processing power of GPUs has increased,
making general-purpose computing on these devices possible.
Additionally, significant improvements in double-precision floating point performance
has broadened the computational opportunities.

In Algorithm 1, we present the pseudo-code for the procedure we used to calculate $s_{i}$
from (2.4). The procedure is adjusted to run efficiently on a GPU, where
a large number of threads execute the same program simultaneously on
different sets of data. This program executed on the device is called the kernel.
Because division is not natively supported by the GPU,
we used inverse multiplication to implement modular reduction.
As we ensure that all operands are less than $2^{51}$, we avoid
code branching and reduce the overall complexity.
For the computations, we also used the fused multiply-add (FMA) operation, natively supported
by the device, where $a\cdot b+c$
is executed as a single instruction. The big advantage of FMA
instructions is that rounding of the product $a\cdot b$ is
performed only once, after the addition. This leaves the intermediate
results of the multiplication in full precision of $104$-bits.

Multiple instances of Algorithm 1 can be executed in parallel, as there are no computational dependencies between
the tasks for each $p$.
In our implementation, kernel execution is limited to $10000$ iterations
after which all operands are reduced modulo $p$. Each process is repeated until $i>\frac{p-3}{2}$,
when all results are returned to the host program and the final
values of $r_{p}$ are calculated. Long kernel execution keeps the arithmetic units occupied, which effectively masks memory latency.
Due to high arithmetic intensity and parallelization capability, we find that Algorithm 1 is suitable for efficient execution in a SIMD environment.

We implemented the algorithm using the OpenCL framework.
As the main computation unit, we used two AMD (ATI) Radeon R9 280x GPUs with a total peak
double-precision performance of $2$ teraFLOPS.
This GPU is based on the Graphics Core Next (GCN) architecture and
has $32$ compute units, each having four $16$-wide SIMD units \cite{AMD}.
In OpenCL terminology, a single kernel instance is called a work-item.
At the hardware level, $64$ work-items are grouped together and executed on an
SIMD unit. This minimum execution unit is called a wavefront.
In a single compute unit, four SIMD units process
instructions from four different wavefronts simultaneously.
Each SIMD unit has an instruction buffer for up to
$10$ wavefronts that can be issued to hide ALU and memory latencies.
Overall, this model of GPU can process up to $81920$ active work-items.
In practice, the number of active wavefronts should be large enough to
hide memory latencies and keep all compute units busy. For the effective utilization, we processed
$65536$ work-items ($1024$ wavefronts) for kernel execution on a single GPU.

Further, we improved the performance by processing two values of $s_{i}$ in the same kernel.
We boosted the performance by an additional 20\% using loop unrolling and instruction pairing.
Using the profiler, we confirmed that the GPU resources were efficiently used.
In this setup, we were able to resolve a block of $262144$ primes in $930$ seconds on average,
for primes in the region of $2^{30}$.
This translates to approximately $1.51\cdot10^{11}$ iterations of the main loop per second, as given in Algorithm 1.
The process took approximately $240$ days to resolve $r_p$ for all $p<2^{34}$.

Comparing this
performance to previous implementations run on a
single i7 class quad-core CPU, our implementation is more than $150$
times faster for $p<2^{32}$. We did not compare similar procedures
on CPUs when $p$ takes values larger than $2^{32}$. We believe that
because of the limitations of the architecture, possible solutions on CPUs might be
even slower.

To verify our results we used two methods. First we selected
random values of $r_{p}$ obtained from the GPU computation,
and verified them on a CPU using a slow but precise procedure
that operates in $128$-bit precision. In this way we verified more than
$200000$ values. It has been reported that in some states, such as when
overclocked, GPUs can generate memory errors \cite{HP}. Thus we also
reran $20$ million randomly selected values of $r_p$
and compared them with the previous results.
For both tests, the results matched exactly.

\begin{algorithm}[H]
\begin{algorithmic}[1]
\Procedure{Kurepa}{$s,i,p$} \State$m\gets (4i^2-2i-2)\bmod p$
\State$k\gets(8i+2)\bmod p$ \State$g\gets1/p$
\State$p_{1}\gets p+1$
\State$c\gets2^{51}+2^{52}$
\State$i_{max}\gets\min\left\{i+10000, \frac{p-3}{2}\right\}$
\While{$i\leq i_{max}$} \State$u\gets s\cdot m$ \State$b\gets
\mbox{FMA}(u,g,c)-c$ \State$s\gets
p_{1}+\mbox{FMA}(s,m,-u)-\mbox{FMA}(p,b,-u)$ \State$m\gets m+k$
\State$k\gets k+8$ \State$i\gets i+1$ \EndWhile \State$s\gets
s\bmod p$ \State\Return $s$ \EndProcedure
\end{algorithmic}
\caption{Computation kernel for $s_i$}
\end{algorithm}

\section{Heuristic considerations}

We strongly believe that there is a counterexample to Kurepa's
conjecture. Heuristic considerations suggest that $!p$ is a random
number modulo $p$ with uniform distribution, so the probability
that $r_p$ has any particular value is approximately $\frac{1}{p}$, and
the sum of reciprocals of the primes diverges. The natural
question then is whether there are infinitely many such counterexamples.
Let us remark that the following considerations are similar to those
in \cite{CDP} and \cite{Z}.

We might expect the number of counterexamples in an interval
$[a,b]$ to be $\sum_{a\leq p\leq b} \frac{1}{p}$. According to
Mertens' second theorem, $\sum_{p\leq x} \frac{1}{p}=\ln \ln
x+c+O(\frac{1}{\ln x})$, where $c\approx 0.261497...$ is the Mertens
constant, and therefore $\sum_{a\leq p\leq b} \frac{1}{p}\approx
\ln(\frac{\ln b}{\ln a})$.

For example, if we consider the interval
$p\in[2^{30},2^{34}]$, then the expected number of primes
$p$ with $\abs{r_p}<l$ is approximately $(2l-1)\ln(\frac{17}{15})$. For
$l=100$, this approach predicts $25$ such primes whereas the
actual number is $23$, and for $l=10000$ it predicts $2503$ such
primes versus the actual number of $2486$.

With such statistical considerations, the probability that there
is no counterexample in an interval $[a,b]$ is $\prod_{a\leq p\leq
b}(1-\frac{1}{p})$. According to Mertens' third theorem
$\prod_{p\leq x}(1-\frac{1}{p})=\frac{e^{-\gamma}}{\ln
x}(1+O(\frac{1}{\ln x}))$, where $\gamma\approx 0.577215...$ is
Euler's constant, and therefore is $\prod_{a\leq p\leq
b}(1-\frac{1}{p})\approx \frac{\ln a}{\ln b}$.

These heuristic considerations predict a $50\%$ chance for a
counterexample in an interval $[2^{t},2^{2t}]$. For example, it
predicts an approximately $20\%$ chance for a counterexample in
the new interval $[1.44\cdot 10^8,2^{34}]$ covered by
this paper. According to this, counterexamples should be rare and
hard to verify. Increasing the upper bound to $p<2^{35}$ gives
only a $3\%$ chance of finding a counterexample, while for a
$50\%$ chance of finding a counterexample we would need to search
$p<2^{68}$, which cannot be attained using existing
hardware and algorithms.

There are many open number theory problems with the same
statistical considerations. For example, searches for
Wieferich primes ($2^{p-1}\equiv 1 \pmod{p^2}$) have yielded only $1093$
and $3511$ for $p<1.4\cdot 10^{17}$ \cite{PG1}, searches for Wilson primes
($(p-1)!\equiv -1\pmod{p^2}$) have yielded only $5$, $13$, and $563$ for
$p<2\cdot 10^{13}$ \cite{CGH}, while no Wall--Sun--Sun primes
($F_{p-(\frac{p}{5})}\equiv 0 \pmod{p^2}$) has been found for
$p<2.8\cdot 10^{16}$ \cite{PG2}. We note that for the listed
problems there are algorithms with better asymptotic time
complexity than $O(n^{2}/\ln n)$, which is the complexity of the best-known
algorithms for searching for a counterexample to Kurepa's conjecture.

\section{Results}

In our search, we calculated and recorded all the residues $r_p$.
It is convenient to use the convention that $r_p$ belongs to
an interval $[-\frac{p-1}{2},\frac{p-1}{2}]$.
In the following table, we list all the values
$\abs{r_p}<100$ for primes
$1.44\cdot 10^{8}<p<2^{34}$.

\vspace{2mm}
\begin{center}
\begin{tabular}{||c|c||c|c||c|c||}
\hline
$p$ & $r_p$ & $p$ & $r_p$ & $p$ & $r_p$ \\
\hline
$145946963$ & $-49$ & $709692847$ & $-38$ & $3730087171$ & $-69$ \\
$171707099$ & $52$ & $758909887$ & $-85$ & $4244621567$ & $58$ \\
$301289203$ & $-57$ & $1023141859$ & $96$ & $4360286653$ & $-48$ \\
$309016481$ & $92$ & $1167637147$ & $67$ & $6855730873$ & $2$ \\
$309303529$ & $26$ & $1250341679$ & $15$ & $8413206301$ & $89$ \\
$345002117$ & $-5$ & $1278568703$ & $-6$ & $9484236149$ & $-64$ \\
$348245083$ & $-83$ & $1283842181$ & $80$ & $11102372713$ & $54$ \\
$353077883$ & $6$ & $1330433659$ & $-75$ & $12024036851$ & $-72$ \\
$441778013$ & $-27$ & $1867557269$ & $-42$ & $14594548571$ & $89$ \\
$473562253$ & $25$ & $2176568417$ & $-60$ & $16541646029$ & $-88$ \\
$499509403$ & $-74$ & $2480960057$ & $-57$ & $16650884357$ & $-82$ \\
$530339209$ & $-24$ & $2649813899$ & $-51$ & $16683898487$ & $91$ \\
$594153589$ & $14$ & $3113484391$ & $-67$ & $$ & $$ \\
$653214853$ & $78$ & $3122927597$ & $-48$ & $$ & $$ \\

\hline
\end{tabular}
\end{center}
\vspace{2mm}

We also considered the generalized left factorial $!^k p=(0!)^k+
(1!)^k+\cdots+((p-1)!)^k$ for small values of $k$. For even $k$, we have
$!^k 3\equiv 1+1+2^k\equiv 0\pmod 3$, and therefore $3\mid !^k
3$. Similarly, $!^k 5\equiv 1+1+2^k+6^k+24^k\equiv
2+2^k+1+(-1)^k\pmod{5}$, and therefore $5\mid !^k 5$ for
$k\equiv 0,3\pmod 4$. Thus, we need to check only $k\equiv 1\pmod
4$.

We searched for the smallest odd prime $p$ such that $p\mid !^k
p$ for all $1<k<100$. The results are listed in the following table,
including a counterexample to Brown's open problem ($k=5$).

\vspace{2mm}
\begin{center}
\begin{tabular}{||r|c||r|c||r|c||r|c||}
\hline
$k$ & \phantom{1234}$p$\phantom{1234} & $k$ & \phantom{1234}$p$\phantom{1234} & $k$ & \phantom{1234}$p$\phantom{1234} & $k$ & \phantom{1234}$p$\phantom{1234}\\
\hline
$5$ & $9632267$ & $29$ & $14293$ & $53$ & $59$ & $77$ & $7852751$\\
$9$ & $29$ & $33$ & $89$ & $57$ & $109$ & $81$ & $229$ \\
$13$ & $97894723$ & $37$ & $29$ & $61$ & $47$  & $85$ & $61$ \\
$17$ & $17203$ & $41$ & $487$ & $65$ & $29$ & $89$ & $104207$\\
$21$ & $27920293$ & $45$ & $233$ & $69$ & $3307$ & $93$ & $29$\\
$25$ & $61$ & $49$ & $859$ & $73$ & 174978647 & $97$ & $18211$\\
\hline
\end{tabular}
\end{center}

\bibliographystyle{amsplain}

\end{document}